\documentclass[twocolumn]{article}
\usepackage[utf8]{inputenc}

\usepackage{amsmath}
\usepackage{amsfonts}
\usepackage{amssymb}
\usepackage{graphicx}
\usepackage{multirow} 

\title{Constructive interpolation points selection in the Loewner framework}
\author{Pierre Vuillemin$^\dagger$ \& Charles Poussot-Vassal$^\dagger$\\
\small 
$^\dagger$ONERA / DTIS, Universit\'e de Toulouse, F-31055 Toulouse, France\\
\small 
\texttt{pierre.vuillemin@onera.fr}, \texttt{charles.poussot-vassal@onera.fr}
}
\date{}

\usepackage{enumitem}
\usepackage{xparse}
\DeclareDocumentCommand \cs {o} {%
  \IfNoValueTF{#1}{\mathbb{C}}{\mathbb{C}^{#1}}%
}
\DeclareDocumentCommand \rs {o} {%
  \IfNoValueTF{#1}{\mathbb{R}}{\mathbb{R}^{#1}}%
}
\DeclareDocumentCommand \ns {o} {%
  \IfNoValueTF{#1}{\mathbb{N}}{\mathbb{N}^{#1}}%
}
\usepackage[ruled,linesnumbered]{algorithm2e}
\usepackage{algorithmic}

\providecommand{\E}[0]{{E}} %
\providecommand{\A}[0]{{A}} %
\providecommand{\B}[0]{{B}} %
\providecommand{\C}[0]{{C}} %
\begin{document}

\maketitle
\begin{abstract}
    This note describes a constructive heuristic to select frequencies of interest within the context of reduced-order modelling by interpolation. The approach is described here through the Loewner framework. Numerical illustrations highlight the benefit it can bring to decrease the required number of interpolation points which is key when the data come from numerically expensive solvers.
\end{abstract}

\section{Introduction}

The Loewner framework (LF) \cite{mayo:framework:2007} is a data-driven method aimed at building a descriptor realisation $\E,\A \in \rs[n\times n], \B \in \rs[n\times p], \C \in \rs[m\times n]$ such that the associated  transfer function $H(z) = C(zE - A)^{-1}B$ interpolates given frequency-domain data $(z_i, \Phi_i)_{i=1}^{N}$, i.e. $H(z_i)=\Phi_i$. When $m,p>1$, the interpolation is tangential. Provided $N$ is large enough, the LF both encodes the minimal McMillan degree $\nu=\mathrm{rank}\,\E$ and the minimal realisation order $n$ to describe the data.

In practical engineering applications, the interpolation points $z_i \in \cs$ are often located on the imaginary axis and the associated data $\Phi_i \in \cs[m \times p]$ then represent the transfer function $G$ of the underlying system, i.e. $\Phi_i = G(z_i)$ - assuming the latter is linear. These data may be evaluated by experiments or a dedicated numerical solver. 

From the authors' experience in (industrial) aeronautic applications (see e.g. \cite{poussot:2020:gla}), the LF has proven to be particularly useful to address the challenges associated with the modelling of complex systems\footnote{Complexity refers here both to the underlying system $G$ and the process to generate the associated data $\Phi_i$ (e.g. through computational fluid dynamics, experiments, etc.) and also to the large dimension $n$ required to catch perfectly all the dynamics they embed.}. 
Indeed, it is simple, numerically efficient and the resulting finite-dimensional linear state-space representation of the interpolant model $H$ is particularly suited to control-engineering applications. Especially as the dimension of $H$ can be decreased even further by a suitable projection in exchange for a loss of accuracy in the interpolation.

Recent work \cite{palitta2021efficient} addresses the issue of memory-management within the LF when the number of data $N$ is large. Yet, considering costly data generated by dedicated high-fidelity solvers, the problematic is reversed as $N$ is likely to be very limited. In that case, the main question lies in the adequate choice of interpolation points $z_i$ to obtain a satisfactory compromise between the accuracy of the model $H$ and the numerical cost to built it.

The LF has been exploited previously in \cite{BeattieCDC:2012} in a fixed-point algorithm where the interpolation points are updated iteratively to fulfil the first-order $\mathcal{H}_2$ optimality conditions. However, the a priori fixed dimension $n$ and the potentially large number of required evaluations of $G$ (in the whole complex plane) before convergence do not suit the considered problem. The AAA algorithm \cite{nakatsukasa:2018:aaa} chooses the next interpolation points based on the worst mismatch with the available data. This supposes that a lot of data are available beforehand.

In this context, we propose in section \ref{sec:constr} a constructive heuristic for the selection of adequate interpolation points. At each step, it exploits the current interpolant model to infer the next frequencies of interest. As highlighted in section \ref{sec:example}, it enables to reduce drastically the required number of function evaluation while reaching an accurate representation of the underlying model.

While writing this note, the  work \cite{cherifi:2021:greedy} came to our attention. The authors develop a similar iterative approach. It differs mainly by the heuristic for selecting the next iteration points which is an interesting alternative.

\section{Constructive Loewner}
\label{sec:constr}
In the sequel, interpolation points $z_i$ are assumed to be chosen solely on the imaginary axis, i.e. $z_i = j \omega_i$, thus allowing to reason only in terms of the frequencies $\omega_i \in \mathbb{R}$. Let us also assume that an interval of interest $\Omega = [\underline{\omega},\overline{\omega}]$ where $\overline{\omega}< \infty$ is available.

The principle of the constructive approach is summarised in algorithm \ref{alg}. Starting from the boundary points of $\Omega$, the set of interpolation points $\mathcal{I}$ is completed at each iteration by the frequency where the current interpolant model $H_k$ has the strongest dynamic. As detailed thereafter, this vague concept can be translated in various ways. But independently, the underlying idea here consists in ensuring that any strong dynamic exhibited by the interpolant model $H_k$ is actually representative of $G$. The algorithm stops when either the maximum allowed number of interpolation points has been reached or the interpolating model does not change more than a prescribed threshold $\epsilon$ from one iteration to the other. A drop in the singular values of the Loewner pencil may also be monitored but as highlighted as it suggests that some data are redundant.

\begin{algorithm}[t]
\begin{algorithmic}[1]
\LinesNumbered
\REQUIRE Initial model $G$, interval $\Omega$, max. number of interpolation points $\overline{r}>2$, tolerance $\epsilon$.
\STATE Let $\mathcal{I} \leftarrow \Omega$
\STATE Evaluate $\Phi_i = G(j\omega_i)$, $\omega_i \in \mathcal{I}$
\STATE $k=0$, $r=\vert \mathcal{I} \vert$
\WHILE{$r \leq \overline{r}$}
\STATE $k \leftarrow k+1$
\STATE Build $H_k$ interpolating $(j\omega_i, \Phi_i)_{i=1}^{r}$
\IF{$\Vert H_{k-1} - H_{k} \Vert / \Vert H_{k-1} \Vert \leq \epsilon$}
\STATE \textbf{break}
\ENDIF
\STATE Find $\omega_{r+1} \in \Omega\setminus \mathcal{I}$ where $H_k(j\omega)$ has the strongest dynamic 
\STATE $\mathcal{I} \leftarrow \mathcal{I} \cup \omega_{r+1}$
\STATE Evaluate $\Phi_{r+1} = G(j\omega_{r+1})$
\STATE $r=\vert \mathcal{I} \vert$
\ENDWHILE
\STATE Return $H_k$
\end{algorithmic}
\caption{Constructive interpolation}
\label{alg}
\end{algorithm}

\paragraph{Notion of strong dynamics. }From a reduction perspective, a low approximation error $\Vert G - H_k \Vert$ between the complex model $G$ and its reduction $H_k$ is, to some extent, related to a look-alike between the singular value plots of both models.

As the error $G - H_k$ is unavailable here, the objective is to make the singular value plot of $H_k$ similar to the one of $G$. The latter being also unknown, the thinking is reversed: for each frequency not yet interpolated $\omega \notin \mathcal{I}$ where $H_k$ exhibits a strong dynamic, one must ensure that $G$ has actually a corresponding strong dynamic.

What a strong dynamic means may largely depend on the application. Flexible structures as those encountered in aeronautic are generally characterised by a resonant frequency response. In that case, the look-alike is obtained by matching both its peaks and valleys. Large static gains, derivative or integral actions are also very common and should also be considered. Therefore, in the sequel, line $9$ of algorithm \ref{alg} will consist in looking first at the largest and lowest gains where $\Vert H'(j\omega)\Vert_2 $ vanishes, then at the largest and lowest slopes.

Many refinement and variations can be imagined but this interpretation has proven to be very effective as highlighted in section \ref{sec:example}. Before that, the implementation details are discussed below.

\paragraph{Practical considerations. }Both the determination of the strong dynamics of $H_k$ and the relative error at line $6$ can be conducted in an exact manner by working with the (small) realisation of the interpolant model. However, we believe that a fully data-driven approach is accurate enough and much simpler to implement.

More specifically, let $\mathcal{W}_{f}$ be a fine discretisation of the interval $\Omega$. Then, the relative error at line $6$ can be replaced by
\begin{equation}
    \tilde{e}_\infty(H_{k-1},H_{k}) =\frac{\max_{\omega \in \mathcal{W}_f}\Vert H_{k-1} (j \omega) - H_k(j\omega) \Vert_2}{\max_{\omega \in \mathcal{W}_f}\Vert H_k(j\omega) \Vert_2}.
    \label{eq:stop_crit}
\end{equation}
Similarly, finding the strong dynamics simplifies as the real domain of search $\Omega \setminus \mathcal{I}$ is reduced to the discrete set $\mathcal{W}_f \setminus \mathcal{I}$. 
The latter can also be refined to avoid clustering of the additional frequencies nearby points already in $\mathcal{I}$. 
The actual research of strong dynamics is done in two consecutive phases. Let $f(\omega) = \Vert H_k(j\omega) \Vert_2$, then,
\begin{itemize}
    \item First, to detect the main peaks and valleys: the derivative $f'(\omega) $, $\omega \in \mathcal{W}_f$ is approximated by finite differences and the frequencies associated with a change of sign are identified. Among those frequencies, the ones with largest or lowest value $f(\omega)$ which are not already in $\mathcal{I}$ are retained.
    \item Should no peak/valley be detected in $\mathcal{W}_f \setminus \mathcal{I}$, then the frequencies such that $ f'(\omega)$ is maximal or minimal are retained.
\end{itemize}
In practice, it is interesting to enrich the interpolation set $\mathcal{I}$ by two points at each iteration. This avoids premature convergence of the approach by stimulating the appearance of new strong dynamics in $H_{k+1}$.



\section{Numerical illustration}
\label{sec:example}
\begin{figure*}
    \centering
    \includegraphics[width=\linewidth,keepaspectratio]{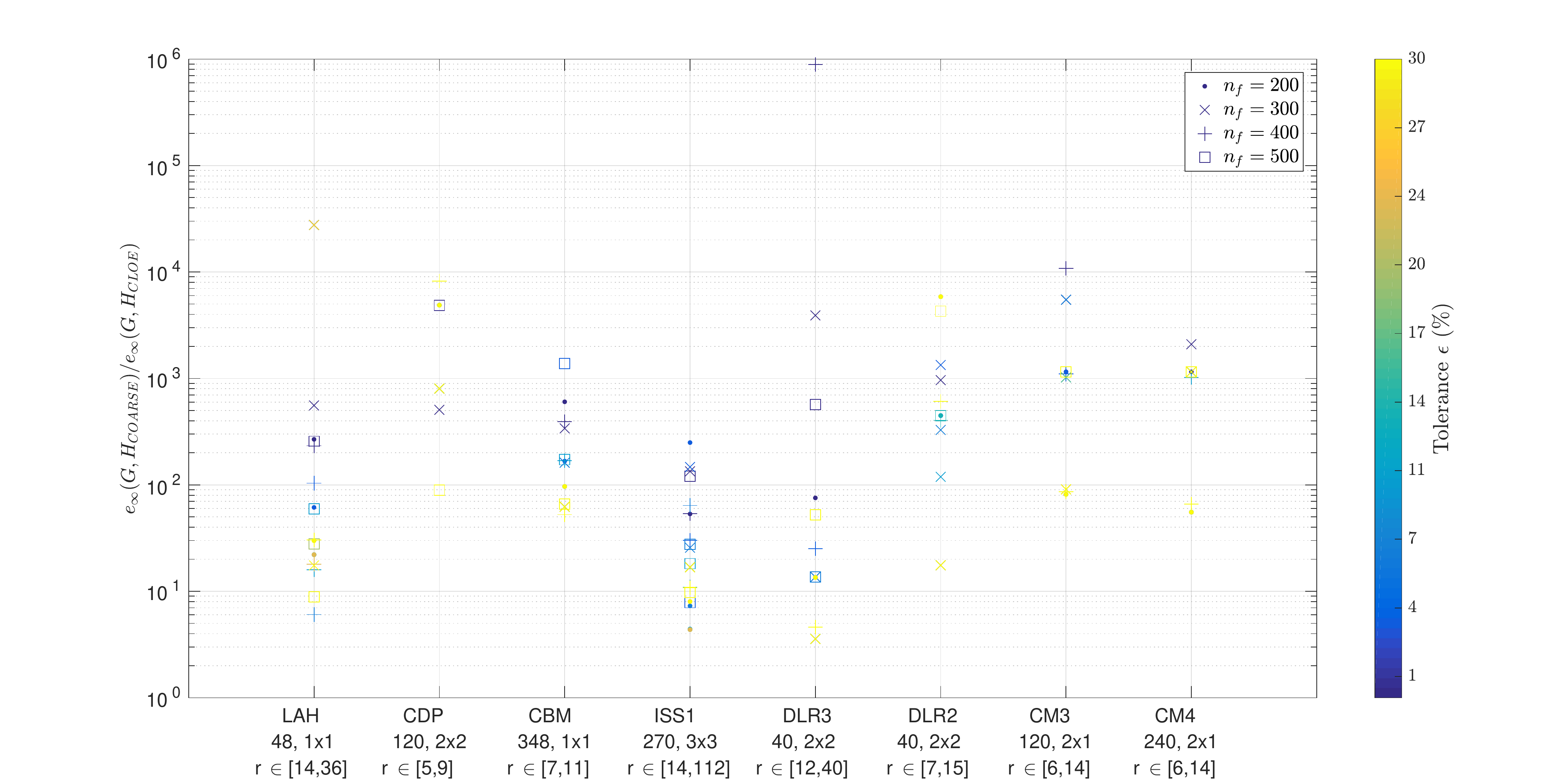}
    \caption{Ratio $e_\infty(G,H_{COARSE})/ e_\infty(G,H_{CLOE})$ for several models which names are reported in abscissa together with their dimension, $m\times p$ and the range of interpolation points $r$ found.}
    \label{fig:cloe1}
\end{figure*}

To highlight the interest of the proposed Constructive Loewner (CLOE), it is applied on several models from Complieb \cite{complib} considering minimum knowledge, i.e. for a common interval $\Omega = [10^{-3}, 10^3]$. The resulting relative $\mathcal{L}_\infty$ error,
\begin{equation}
    e_\infty(G,H) = \Vert G - H \Vert_\infty / \Vert G \Vert_\infty,
\end{equation}
is compared to the one obtained with the Loewner model (referred as coarse in the sequel) built with the same number of points but spread logarithmically over $\Omega$. Within CLOE, the set $\Omega$ is discertised logarithmically with a number of samples $n_f$ varying from $200$ to $500$ and the tolerance $\epsilon$ varies from $1\%$ to $30\%$. The results are reported in figure \ref{fig:cloe1}.

\begin{figure}
    \centering
    \includegraphics[width=\linewidth,keepaspectratio]{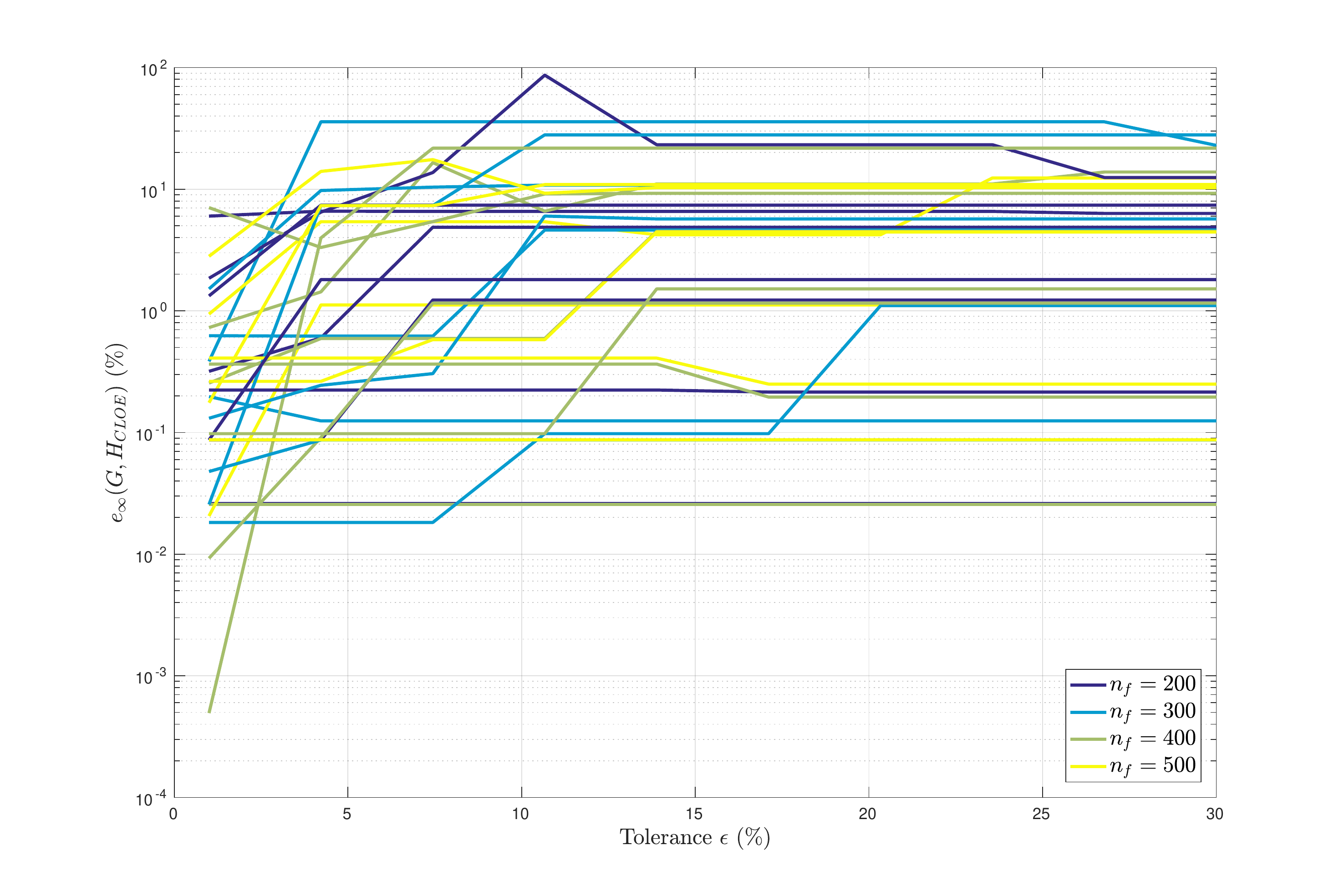}
    \caption{Relative error of $H_{CLOE}$ against the tolerance $\epsilon$ for various values of $n_f$ and all the considered test models.}
    \label{fig:err_mdls}
\end{figure}

One can notice that the proposed heuristic is an improvement in comparison to a logarithmic selection of the interpolation points as the resulting error is lower in all the considered cases\footnote{Note that several points are overlapping as the number of interpolation points found by CLOE is evolving by steps.}. This is partly luck as those performances are not to be expected for the largest values of $\epsilon$ (see the iterations on the LAH model thereafter). No clear trend appears when the dimension $n_f = \vert \mathcal{W}_f \vert$ or the tolerance $\epsilon$ vary.
However, when looking at the relative error of the reduced-order model in figure \ref{fig:err_mdls}, one can see that a lower approximation error is reached when the tolerance $\epsilon$ is chosen below $5 \%$. As illustrated thereafter, it is indeed important to stop the algorithm when the interpolant model does not evolve anymore throughout the iterations.

\begin{figure*}
\centering 

\includegraphics[width=0.49\linewidth, keepaspectratio]{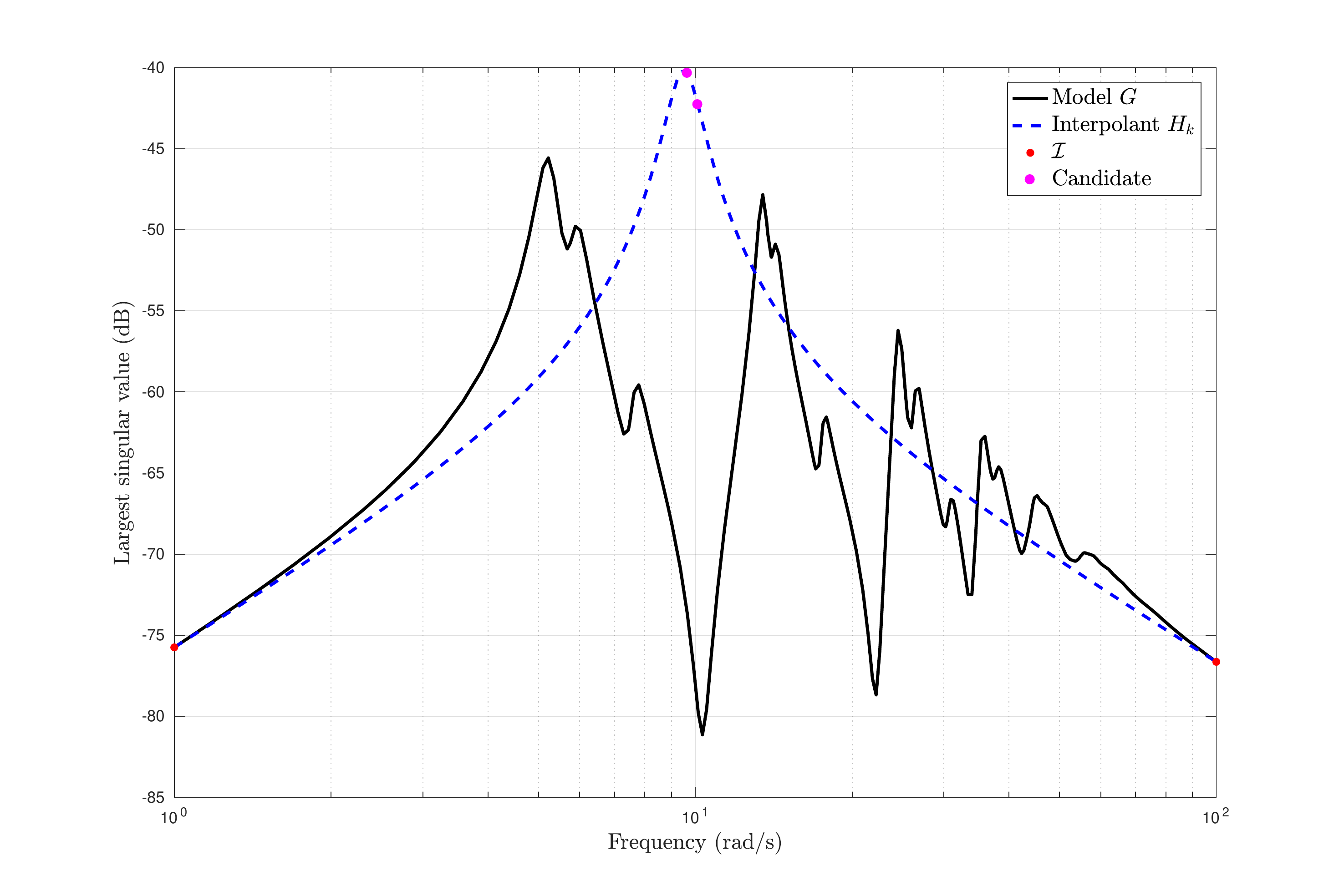}
\includegraphics[width=0.49\linewidth, keepaspectratio]{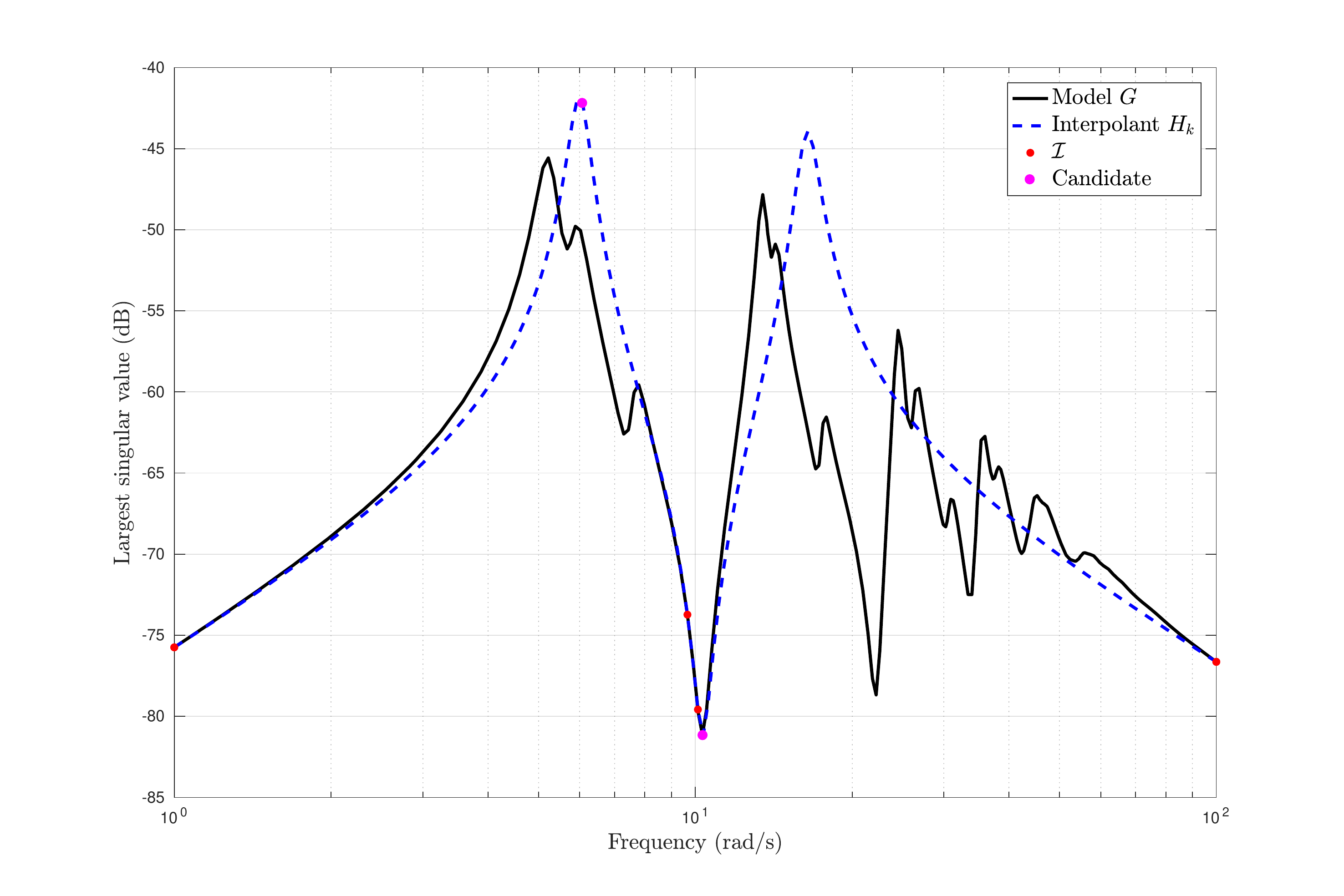}

\includegraphics[width=0.49\linewidth, keepaspectratio]{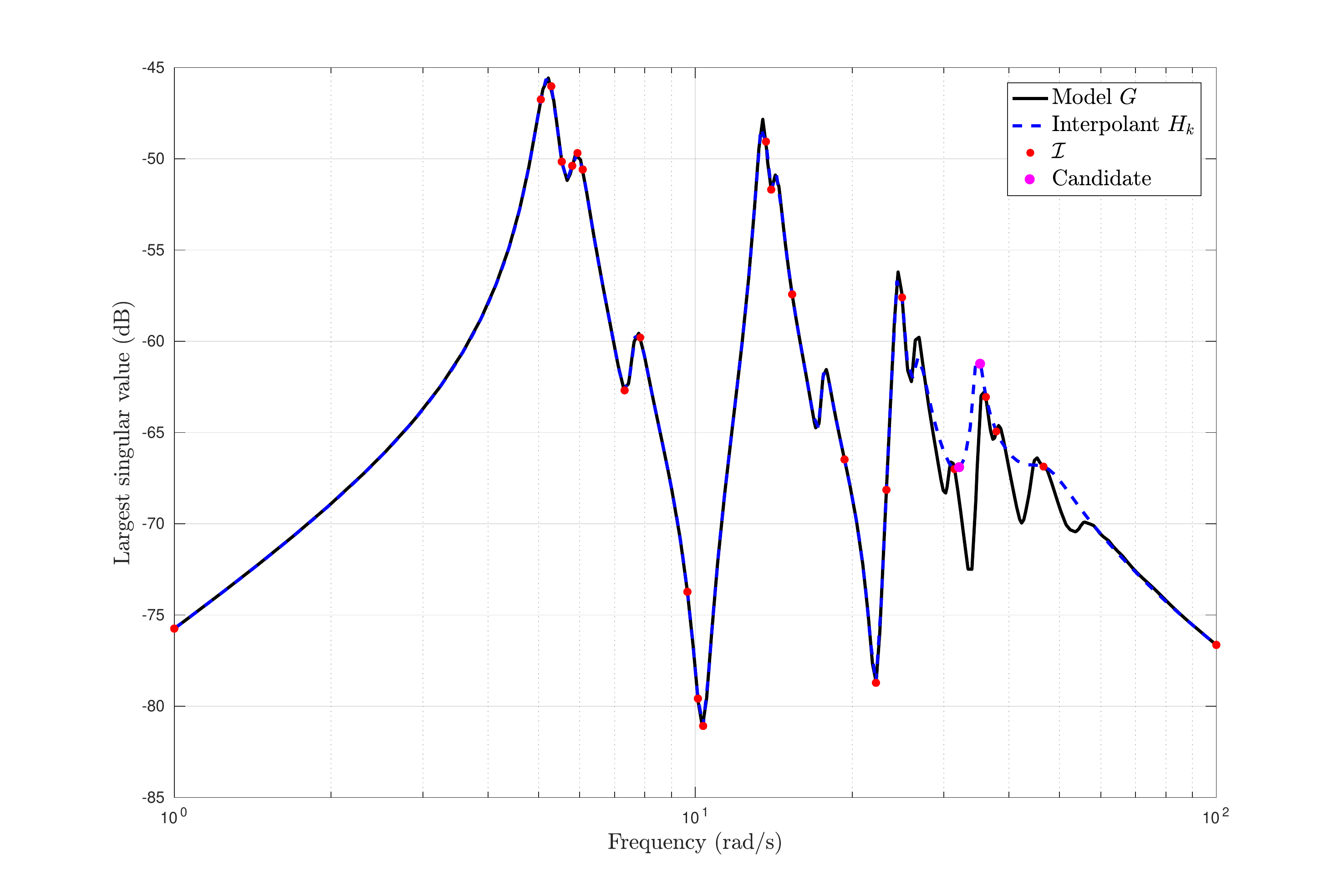}
\includegraphics[width=0.49\linewidth, keepaspectratio]{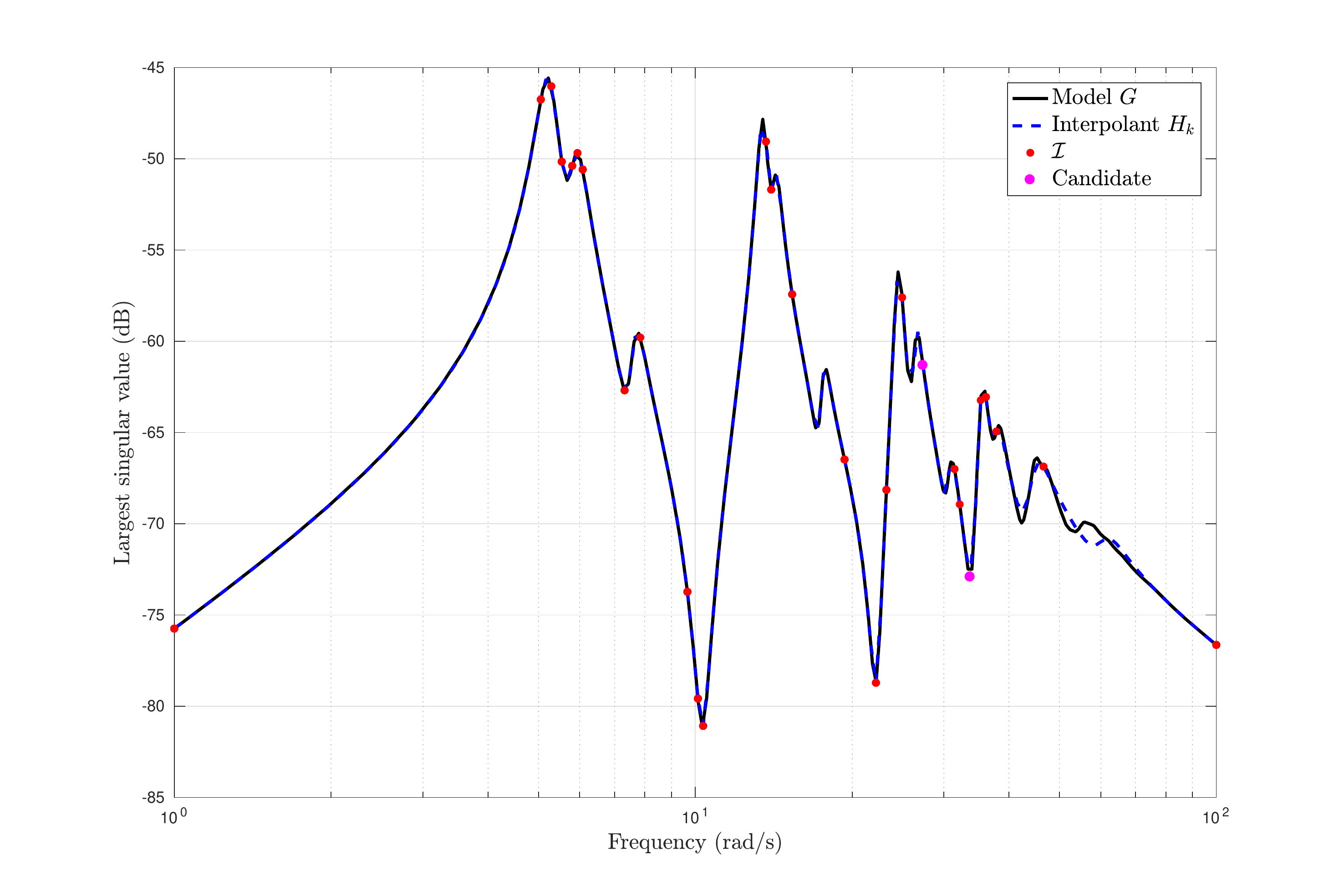}

\caption{Frequency responses of the LAH model (black) together with the interpolant model $H_k$  obtained with CLOE at iteration $k$ (dashed blue). The first two iterations are displayed at the top and the last two iterations at the bottom. Selected interpolation points (set $\mathcal{I}$) are in red and candidate points are in magenta. The parameter $\epsilon$ is set to $5\%$.}
\label{fig:bode}
\end{figure*}
 
\begin{figure*}
    \centering
    \includegraphics[width=\linewidth,keepaspectratio]{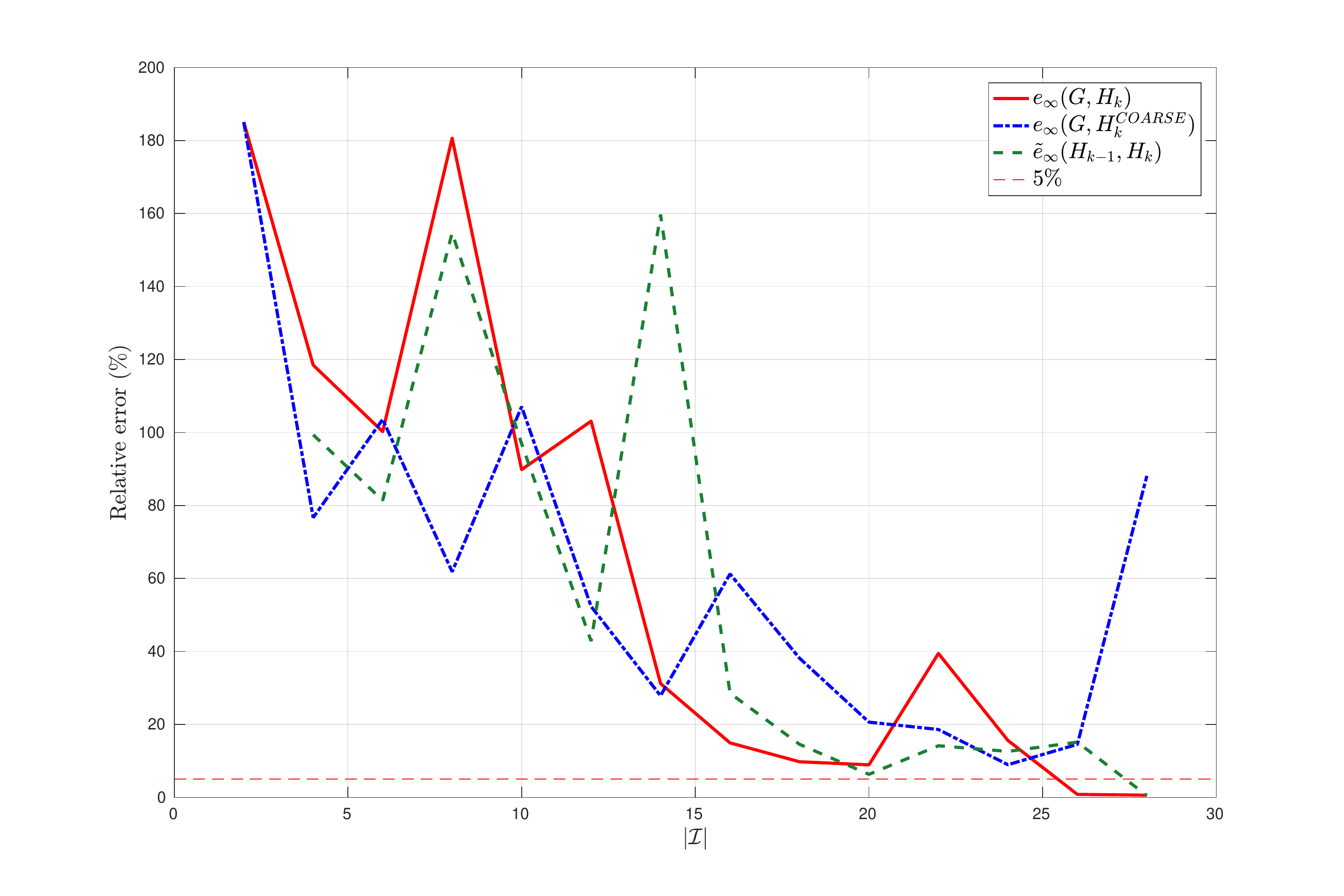}
    \caption{Relative errors between the LAH model $G$ and the models $H_k$ and $H_k^{COARSE}$ obtained respectively with CLOE at iteration $k$ or with the coarse grid with the same number of points. The stopping criterion \eqref{eq:stop_crit} is also reported. }
    \label{fig:err}
\end{figure*}

To illustrate the evolution of the model within CLOE, the interpolating model $H_k$ is shown together with the initial model $G$ in figure \ref{fig:bode} at various iterations. The evolution of the relative approximation error throughout the iterations is reported in figure \ref{fig:err} together with the error of corresponding coarse model and the value of the stopping criterion \eqref{eq:stop_crit}.

In figure \ref{fig:bode} top left, the initial model is interpolated at the boundary of the frequency interval, and the candidate interpolation points are near the only resonance of $H_2$: one candidate point is the detected peak and the other one the most negative slope. During the next iteration, a peak and a valley are detected. The two last iterations show that the main peaks are reproduced by the interpolant model and that smaller dynamics are progressively being matched. 

One can see in figure \ref{fig:err} that the approximation errors do not decrease monotonically as the number of interpolation points increases. Still, the stopping criteria $\tilde{e}_\infty$ appears to be relevant as it follows the trend of the real approximation error reached by CLOE. The importance of the tolerance $\epsilon$ is highlighted here by the fact that throughout the iterations, the error can be larger than with the coarse grid, especially when $\tilde{e}_\infty(H_{k-1},H_k)$ is large.

\section{Conclusion}

This note details an heuristic approach for constructive Loewner interpolation. The number of interpolation points in gradually increased in order to try to improve the quality of the interpolant model while limiting the number of required evaluation of the underlying system $G$. The choice of new interpolation points is the key of the approach. It relies solely on the interpolant model and aims at ensuring that any strong dynamic it contains is actually representative of $G$. These dynamics are identified in the frequency domain by looking for the peaks, valleys and large slope areas of the singular values plot.

In spite of its simplicity, this heuristic has shown to be effective in order to determine an accurate model while limiting the number of required data from the underlying system. It opens interesting perspectives for the modelling from parsimonious frequency-domain data provided the latter can be selected.

\bibliographystyle{plain} 
\bibliography{piblio}
\end{document}